\documentclass[10pt,a4]{amsart}
\usepackage{amsmath,amssymb}
\makeatletter
\renewcommand{\section}{\@startsection{section}{1}{0pt}{20pt}{6pt}{\large\bfseries}}
\makeatother

\usepackage{enumerate}

\numberwithin{equation}{section}

\theoremstyle{plain}
  \newtheorem{thm}{Theorem}[section]

 \theoremstyle{definition}

\newcommand{\R}{\mathbb{R}}

\renewcommand{\Re}{{\mathfrak{Re}}}
\newcommand{\C}{\mathbb{C}}
\newcommand{\Q}{\mathbb{Q}}
\newcommand{\E}{\mathbb{E}}
\newcommand{\e}{{\bf{e}}}
\newcommand{\tp}{{\overline{\psi}}}
\renewcommand{\P}{{\mathbb{P}}}

\renewcommand{\O}{\mathcal{O}}

\newcommand{\I}{\mathcal{I}}
\newcommand{\Ip}{\mathcal{I}_{\psi}}

\newcommand{\Ig}{\I_{\psi_{\gamma}}}
\newcommand{\Og}{\O_{\psi_{\gamma}}}

\newcommand{\Id}[1]{{{\mathbb{I}}}_{\{#1\}}}

\newcommand{\asa}[1]{\quad \textrm{as } \: #1\rightarrow \infty}

\newcommand{\ang}{a_n(\psi_{\gamma})}

\begin{document}
\bibliographystyle{plain}
\title{{\small{Probability Theory}}  \\ \vspace{1cm}Law of the exponential functional of one-sided L\'evy processes and Asian options}
\author{Pierre Patie}
\address{Institute of Mathematical Statistics and Actuarial Science, University of Bern\\
Alpeneggstrasse, 22, CH-3012 Bern, Switzerland.}
\email{patie@stat.unibe.ch}

\begin{abstract}
The purpose of this note is to describe, in terms of a power series,  the distribution function of the exponential functional, taken at some  independent exponential time,  of a  spectrally negative L\'evy process $\xi=(\xi_t,t\geq0)$ with unbounded variation.  We also derive a Geman-Yor type formula for Asian options prices in a financial market driven by $e^\xi$.

 \vspace{0.5cm}

\noindent \textsc{R\'esum\'e.} \textbf{Loi de la fonctionnelle exponentielle de processus de L\'evy  assym\'etriques et options asiatiques}. L'object de cette note est de donner une repr\'esentation, en terme d'une s\'erie enti\`ere, de la distribution de la fonctionnelle exponentielle, consid\'er\'ee en un temps exponentiel ind\'ependent,  d'un processus de L\'evy $\xi$ spectrallement n\'egatif,  \`a variation infinie et pouvant \^etre tu\'e. Nous en d\'eduisons une formule du type Geman-Yor pour le prix des options asiatiques dans un march\'e financier dirig\'e par $e^\xi$.
 \end{abstract}

\maketitle

\section*{Version fran\c{c}aise abr\'eg\'ee} Soit $\xi=(\xi_t)_{t\geq0}$ un processus de L\'evy \`a valeurs r\'eelles  spectralement n\'egatif et dont les trajectoires sont \`a variation  infinie. Cela signifie que $\xi$ est un processus dont les accroissements sont stationnaires et ind\'ependants et par ailleurs le processus n'effectue que des sauts n\'egatifs. Il est bien connu que la loi du processus $\xi$ est caract\'eris\'ee par la loi de la variable al\'eatoire $\xi_1$ et par cons\'equent par l'exposant de Laplace de cette derni\`ere que nous \'ecrivons $\psi$. Sous les conditions {\bf{H}}, donn\'ees dans le corps de la note, nous proposons de d\'ecrire la loi de la variable al\'eatoire
$ \Sigma_{\e_q} =\int_0^{\e_q}e^{\xi_s}ds,$
o\`u $\e_q$ est une variable al\'eatoire ind\'ependente de $\xi$ et suivant une loi exponentielle de param\`etre $q\geq0$, o\`u nous comprenons que $\e_0=\infty$. La distribution  de $\Sigma_{\e_q}$ appara\^it dans diff\'erents domaines des probabilit\'es et \'egalement dans diff\'erents champs des math\'ematiques appliqu\'ees. Malheureusement, la connaissance explicite de cette loi se r\'eduit \`a quelques cas particuliers, dont celui du mouvement brownien avec d\'erive. A cet effet, nous indiquons  l'excellent papier de Bertoin et Yor \cite{Bertoin-Yor-05} o\`u une description de ces cas particuliers et des enjeux sous-jacents \`a l'\'etude de la loi de $\Sigma_{\e_q}$ sont d\'etaill\'es.   Le  r\'esultat principal que nous \'enon\c{c}ons dans cette note consiste en une  repr\'esentation  de la loi de $\Sigma_{\e_q}$ en terme d'une s\'erie enti\`ere dont les coefficients sont d\'efinis \`a l'aide de l'exposant de Laplace $\psi$.
Une cons\'equence int\'eressante de cette repr\'esentation est l'obtention d'une formule  pour le prix des options asiatiques dans un march\'e financier dirig\'e par $e^{\xi}$. Ce r\'esultat g\'en\'eralise la formule de Geman et Yor \cite{Geman-Yor-92} obtenue dans le cadre du mod\`ele de Black-Scholes.

\vspace{5mm}
\section{Introduction}
Let $\xi=(\xi_t)_{t\geq0}$  be a real-valued spectrally negative L\'evy process with unbounded variation and we denote its law by $\P_y$ $(\P=\P_0)$ when $\xi_0=y\in \R$.
That means that $\xi$ is a process with stationary and independent increments having only negative jumps and its right continuous paths with left-limits are of infinite variation on every compact time interval a.s. We refer to the excellent monographs \cite{Bertoin-96} and \cite{Sato-99} for background. It is well-known that the law of $\xi$ is  characterized by its one dimensional distributions and thus by the Laplace exponent $\psi : \R^+ \rightarrow \R$ of the random variable $\xi_1$ which admits the following   L\'evy-Khintchine representation
 \begin{eqnarray*}
\psi(u) = bu + \frac{\sigma}{2} u^2 + \int_{-\infty}^0 (e^{u r} -1
-ur\Id{|r|<1} )\nu(dr),\quad u\geq0,
\end{eqnarray*}
where $ b\in \R, \:\sigma \geq 0$ and the measure
$\nu$ satisfies the integrability condition $\int^{0}_{-\infty}(1
\wedge r^2 )\:\nu(dr) <+  \infty$. Since we excluded the case when $\xi$ has finite variation, the condition $\int^{0}_{-\infty}(1
\wedge r)\:\nu(dr) <+  \infty$ is not allowed.  Note, see e.g.~\cite[VII, Corollary 5]{Bertoin-96}, that the property that $\xi$ has unbounded variation is equivalent to   the following  asymptotic
\begin{equation} \label{eq:ll}
\lim_{u\rightarrow \infty} \frac{\psi(u)}{u}=+\infty.
\end{equation}
The aim of this note is to describe the distribution function of the so-called exponential functional
\[ \Sigma_{\e_q} =\int_0^{\e_q}e^{\xi_s}ds\]
where $\e_q$ is a random variable, independent of $\xi$, which is exponentially distributed with parameter $q\geq 0$, where we understand $\e_0=\infty$. In particular, if $q=0$, the strong law of large numbers for L\'evy processes gives the following equivalence
\[ \Sigma_{\e_0}<\infty \quad \textrm{ a.s.}\quad \Longleftrightarrow \quad  \E[\xi_1]<0\]
and we refer to the  paper of Bertoin and Yor \cite{Bertoin-Yor-05} for alternative conditions, references on the topic and for motivations for studying the  law of $\Sigma_{\e_q}$. We simply mention that this positive random variable appears in various
fields such as diffusion processes in random environments,
fragmentation and coalescence processes, the classical moment
problems, mathematical finance and astrophysics. We also point out that the law of $\Sigma_{\e_q}$ was known only for the Brownian motion with drift, see Yor's monograph \cite{Yor-01}, and a  few other isolated cases, see e.g.~Carmona et al.~\cite{Carmona-Petit-Yor-97}, Gjessing and Paulsen \cite{Gjessing-Paulsen-97} and Patie \cite{Patie-06-poch}.  We are now ready to summarize the conditions which will be in force throughout the remaining part of this note.

\vspace{0.2cm}
\noindent {\bf{H}}: \eqref{eq:ll} holds and either $q>0$ or $q=0$ and $\E[\xi_1]<0$.
\vspace{0.2cm}

\noindent We mention that in \cite{Patie-abs-08} the case when the condition \eqref{eq:ll} does not hold is also considered. The remaining part of this Note  is organized as follows. In the next Section, we state the representation of the distribution of $\Sigma_{\e_q}$ in terms of a power series. In Section 3, we derive a Geman-Yor type formula for the price of Asian options in a spectrally negative L\'evy market. We end up this note by revisiting the Brownian motion with drift case.

\section{Main result}
Let us start by  recalling  some basic properties of the Laplace exponent $\psi$, which can be found in \cite{Bertoin-96}. First, it is plain that $\lim_{u \rightarrow
\infty}\psi(u)=+\infty$ and $\psi$ is convex. Note that $0$ is
always a root of the equation $\psi(u)=0$. However, in the case
$\E[\xi_1] <0$, this equation admits another positive root, which we denote
by $\theta$. Moreover, for any $\E[\xi_1]  \in [-\infty,\infty)$, the function $u\mapsto \psi(u)$ is continuous and increasing on $[\max(\theta,0),\infty)$. Thus, it  has a
well-defined inverse function $\phi:[0,\infty)\rightarrow
[\max(\theta,0),\infty)$ which is also continuous and  increasing.
Finally, we write, for any $u\geq0$ and $q >0$, $\tp(u)=\psi(u)-q$,  and set the following notation
$$
\psi_{\theta}(u)=\psi(u+\theta) \quad \textrm{ and }\quad
\tp_{\phi(q)}(u)=\tp(u+\phi(q)).
$$
Recalling that $\psi(\theta)=0$ and observing that $\tp(\phi(q))=0$, the mappings $\psi_{\theta}, \tp_{\phi(q)}$ are plainly Laplace  exponents of conservative L\'evy processes. We also point out that $\psi_{\theta}'(0^+)=\psi'(\theta^+)>0$ and $\tp_{\phi(q)}'(0^+)=\psi'(\phi(q))=\frac{1}{\phi'(q)}>0$. In order to present our result in a compact form,  we write
$$ 
\gamma=
\begin{cases}
\phi(q) & \quad \textrm{if } q>0, \\
\theta &  \quad  \textrm{otherwise,}
\end{cases} \qquad
\textrm{ and} \qquad
\psi_{\gamma}=
\begin{cases}
\tp_{\phi(q)} & \quad \textrm{if } q>0, \\
\psi_{\theta}&  \quad  \textrm{otherwise.}
\end{cases}
$$ 
Next, set $a_0=1$ and $\ang=\left( \prod_{k=1}^n\psi_{\gamma}( k)\right)^{-1},\: n=1,2,\ldots$
In \cite{Patie-06c}, the author introduced the following power series
\begin{equation*} \Ig(z)= \sum_{n=0}^{\infty}
 \ang z^{n}
\end{equation*}
and showed by means of classical criteria that the mapping $z\mapsto \Ip(z)$ is an entire function. We refer to \cite{Patie-06c} for interesting analytical properties enjoyed by these power series and also for  connections with well known special functions, such as, for instance, the modified Bessel functions, confluent hypergeometric functions and several generalizations of the Mittag-Leffler functions. To simplify the notation, we introduce the following definition, for any $z\in \C$,
\begin{eqnarray*}
\Og( z) &=& \Ig(e^{i\pi}z).
\end{eqnarray*}
 Next, let $G_{\kappa}$ be a Gamma random variable independent of $\xi$, with parameter $\kappa>0$. Its density is given by $
g(dy) = \frac{e^{-y}y^{\kappa-1}}{\Gamma(\kappa)} dy, y>0$, with $\Gamma$ the Euler gamma function. Then, in \cite{Patie-08a},  the author suggested the following generalization
\begin{eqnarray*}
\Og(\kappa ;z) = \E\left[\Og\left(G_{\kappa}z\right)\right].
\end{eqnarray*}
Thus,  by means of the integral representation of the Gamma function $\Gamma(\rho)=\int_0^{\infty}e^{-s}s^{\rho-1}ds,\: \Re(\rho)>0,$ and an argument of dominated convergence, we obtain the following power series representation
\begin{eqnarray} \label{eq:defo}
\Og(\kappa ;z) = \frac{1}{\Gamma(\kappa)}\sum_{n=0}^{\infty}
 (-1)^n \ang
 \Gamma(\kappa+n) z^{n}
\end{eqnarray}
which is easily seen to be, under the condition \eqref{eq:ll}, an entire function in $z$.  Note also, by the principle of analytical continuation, that the mapping $\rho \mapsto \Og(\rho ;z)$ is an entire function for arbitrary $z\in \C$. We are now ready to state the main result of this note.
\begin{thm} \label{thm:2}
Assume that {\bf{H}} holds and  write $S(t)=\P(\Sigma_{\e_q}\geq t),t>0.$
Then, there exists a constant $C_{\gamma}>0$ such that
\begin{eqnarray*}
\O_{\psi_{\gamma}}( \gamma;t)&\sim& \frac{t^{-  \gamma}}{C_{\gamma}}   \asa{t},
\end{eqnarray*}
and one has, for any $t>0$,
\begin{eqnarray*}
S(t)&=& C_{\gamma}t^{-  \gamma} \O_{\psi_{\gamma}}( \gamma;t^{-1}). \label{eq:pd}
\end{eqnarray*}
Moreover, the law of $\Sigma_{\e_q}$ is absolutely continuous with a density, denoted by $s$, given by
\begin{eqnarray*}
s(t)&=& \gamma C_{\gamma}  t^{-  \gamma-1} \O_{\psi_{\gamma}}(1+ \gamma;t^{-1}), \: t>0.
\end{eqnarray*}
\end{thm}
We now sketch the main steps used for proving the Theorem and further details will be provided in \cite{Patie-abs-08}.
To  this end, we denote by  $X=((X_t)_{t\geq0}, (\Q_x)_{x>0})$ a  $1$-self-similar Hunt process with values in $ [0,\infty)$. It means that $X$ is a right-continuous strong Markov process with quasi-left continuous trajectories and  $X$ enjoys the following self-similarity property: for each $c>0$ and $x\geq0$,
\[  \textrm{ the law of the process } (c^{-1}X_{ct})_{t\geq0}, \textrm{under } \Q_x, \textrm{ is  } \Q_{x/c}.\]
 This class of processes was introduced and studied by Lamperti \cite{Lamperti-72}. In particular, Lamperti proved that there is a bijective correspondence
between $[0,\infty)$-valued self-similar Markov processes and real-valued L\'evy
processes.  Moreover, we deduce from one of the Lamperti zero-one laws that under the condition $\bf{H}$,   $\Q_x(T_0<\infty)=1$,  for all $x>0$, where
\[T_0=\inf\{s>0; X_s=0\}\]
is the absorption time of $X$.
Then, the first key step in our proof is the following identity
\begin{eqnarray*}
\Q_{xt^{-1}}(H_0<\infty)  &=&\Q_x( T_0\geq t), \quad x,t>0,
 \end{eqnarray*}
where
\begin{equation} \label{eq:hou}
H_0=\inf\{s>0; \: U_s:=e^{s}X_{(1-e^{- s})}=0\}.
\end{equation}
That is $H_0$ is the absorption time of $U,$ the  so-called Ornstein-Uhlenbeck process associated to $X$ of parameter $-1$. It is also a transient Hunt process on $[0,\infty)$. Note that the identity above is easily obtained by means of the self-similarity property of $X$ and a simple time change. Observing that for $t=1$, the mapping $x\mapsto \Q_{x}(H_0<\infty)$, defined on $\R^+$, is an increasing invariant function for the semigroup of $U$, we have thus transformed a parabolic integro-differential problem into an elliptical one.  Then, specializing on the case when $X$ is  the self-similar Markov process associated to $\xi$ via the Lamperti bijection, we adapt, to the current situation, some devices which have been recently developed by the author in \cite{Patie-06c} and \cite{Patie-08a}, for describing the invariant function of stationary Ornstein-Uhlenbeck processes. However, several  issues arise when dealing with the transient ones. Indeed, in the stationary case, it is an easy matter to derive some basic but essential properties, such as positivity and monotonicity, of the invariant functions as they are  expressed in terms of analytical power series having only positive coefficients. As one observes from \eqref{eq:defo}, it is not as straightforward in the transient case and, for instance, some information about the location and the monotonicity of the real zeros, with respect to the argument $\kappa$, of this entire function is required.  This is achieved by combining probabilistic techniques with basic tools borrowed from complex analysis.

\section{A Geman-Yor type formula for Asian options}
In \cite{Geman-Yor-92}, Geman and Yor derived the price of the so-called Asian option in the Black-Scholes market, i.e.~when the dynamics of the asset price is given as the exponential of a Brownian motion with drift.  More specifically, they compute,  for any $K>0$ and $y\in \R$, the following functional
\[ A(y,K,q)= \E_y[(\Sigma_{\e_q}-K)^+]  \]
in terms of the confluent hypergeometric function, in the case $\xi$ is a Brownian motion with drift.
Before stating the generalization of  Geman-Yor formula, let us point out the following identity $ A(y,K,q)= e^{y}A(0,Ke^{-y},q)$, which follows readily from the translation invariance of L\'evy processes. We also mention that the fundamental theorem of asset pricing, see Delbaen and Schachermayer \cite{Delbaen-Schachermayer-94} requires that the discounted value of the asset price is a (local)-martingale under an equivalent martingale measure. However, it is easily checked that  with the condition $\psi(1)=r$, where $r>0$ is the risk-free rate, one may  carry out  the pricing under $\P$. We are now ready to state the following.
\begin{thm} With the notation used in Theorem \ref{thm:2}, for any $K>0$ and $q>\psi(1)$, we have
\begin{eqnarray*}
\E[(\Sigma_{\e_q}-K)^+]&=&  \frac{C_{\phi(q)}}{(\phi(q)-1)} K^{1-  \phi(q)} \O_{\psi_{\phi(q)}}(\phi(q)-1;K^{-1}).
\end{eqnarray*}
\end{thm}
It should be  emphasized  that, in general, the formula above is not obtained, from Theorem \ref{thm:2}, by a simple term by term integration. The details of the proof of this result are provided  in \cite{Patie-asian-08}.

\section{The Brownian motion with drift case revisited} \label{sub:bess}
We consider $\xi$ to be a $2$-scaled Brownian motion with drift $2b \in \R$ and killed at some independent exponential time of parameter $q>0$, i.e.~$\tp(u)=2u^2+2b u-q$ and $2\phi(q)=\sqrt{2q+b^2}-b$. Note that $\tp_{\phi(q)}(u)=2u^2+(2b+\phi(q)) u$. Its associated 
self-similar process $X$ is well known to be a  Bessel process of index
$b$ killed at a rate $q\int_0^{t}X_s^{-2}ds$. Moreover, we obtain, setting $\varrho=b+2\phi(q)$,
\begin{eqnarray*}
\O_{\tp_{\phi(q)}}\left(\rho;x\right)&=&\frac{\Gamma(\varrho+1)}{\Gamma(\rho)}\sum_{n=0}^{\infty} (-1)^n\frac{\Gamma(\rho+n)}{n!\Gamma(n+\varrho+1)}(x/2)^n\\
&=&\Phi\left(\rho, \varrho+1;-x/2\right)
\end{eqnarray*}
where $\Phi$ stands for the confluent hypergeometric function. We refer to Lebedev \cite[Section 9]{Lebedev-72} for useful properties of this function. Next, using the following  asymptotic \begin{eqnarray*}
\Phi\left(\rho, \varrho+1;-x\right) \sim \frac{\Gamma(\varrho+1)}{\Gamma(\varrho+1-\rho)} x^{-\rho} \asa{x},
\end{eqnarray*}
we get that $C_{\phi(q)}=\frac{\Gamma(\varrho+1-\phi(q))}{2^{\phi(q)}\Gamma(\varrho+1)}$.
 Thus, we obtain, recalling that, for any $q>0$, $\varrho-\phi(q)=b+\phi(q)>0$,
\begin{eqnarray*}
s_{\phi(q)}(t)&=&\phi(q) \frac{\Gamma(\varrho+1-\phi(q))}{2^{\phi(q)}\Gamma(\varrho+1)}t^{-\phi(q)-1}\Phi\left(1+\phi(q), \varrho+1;-(2t)^{-1}\right)\\
&=&\frac{\varrho-\phi(q)}{2^{\phi(q)}\Gamma(\phi(q))}t^{-\phi(q)-1} \int_0^{1}e^{-\frac{u}{2t}}(1-u)^{\varrho-\phi(q)-1}u^{\phi(q)}du
\end{eqnarray*}
which is the expression  \cite[(5.a) p.105]{Yor-01}.

We end up this note by  mentioning that in \cite{Patie-abs-08} some further known and new examples are detailed. In particular, we recover the recent result obtained by  Bernyk et al.~\cite{Bernyk-Dalang-Peskir-08} regarding the density of the maximum of regular spectrally positive stable processes.

\end{document}